\documentclass{amsart}
\usepackage[utf8]{inputenc}

\usepackage{amsmath}
\usepackage{amssymb}
\usepackage{enumitem}
\usepackage{amsthm}
\usepackage[dvipsnames]{xcolor}
\usepackage[pagebackref,hypertexnames=false, colorlinks, citecolor=BrickRed,linkcolor=MidnightBlue, urlcolor=BrickRed]{hyperref}

\usepackage{dsfont}
\usepackage{mathtools}

\usepackage{etoolbox}
\apptocmd{\sloppy}{\hbadness 10000\relax}{}{}
\apptocmd{\sloppy}{\vbadness 10000\relax}{}{}

\usepackage{orcidlink}

\usepackage[pagewise]{lineno}

\newtheorem{theorem}{Theorem}[section]
\newtheorem{lemma}[theorem]{Lemma}
\newtheorem{corollary}[theorem]{Corollary}
\theoremstyle{definition}
\newtheorem{definition}[theorem]{Definition}
\newtheorem{remark}[theorem]{Remark}
\newtheorem{example}[theorem]{Example}
\newtheorem{observation}[theorem]{Observation}
\newtheorem{problem}[theorem]{Problem}

\newcommand{\unitization}{\mathds 1}
\renewcommand{\restriction}{\mathbin{\upharpoonright}}

\title{Schoenberg correspondence for multifaced independence}
\author{Malte Gerhold}

\address{Greifswald University, Institute of Mathematics and Computer Science, Greifswald, Germany}
\thanks{
This work was supported by the German Research Foundation (DFG) grant nos.\ 465189426 and 397960675; it was carried out as a postdoctoral researcher at Saarland University, during the tenure of an ERCIM ``Alain Bensoussan'' Fellowship Programme at NTNU Trondheim, and as a postdoctoral scientific employee at University of Greifswald.}
\urladdr{\href{https://sites.google.com/view/malte-gerhold/}{\url{https://sites.google.com/view/malte-gerhold/}}\newline {\textnormal{{\it ORCID}:}} {\orcidlink{}}\,\url{https://orcid.org/0000-0003-4029-1108}}

\begin{document}


\begin{abstract}
    We extend the Schoenberg correspondence for universal independences by Schürmann \& Voß to the multivariate setting of Manzel \& Schürmann, covering, e.g., Voiculescu's bifreeness as well as Bo{\.z}ejko \& Speicher's c-free independence. At the same time, we free the proof in the univariate situation from its dependence on Muraki's ``5 Independences Theorem''.
\end{abstract}
  
\maketitle

\vspace{-1em}

\section{Introduction}

The term \emph{Schoenberg correspondence} refers to theorems that characterize one-parameter semigroups which are (in some sense) positive by a corresponding condition on their generators. Schoenberg's original version (in today's language this is the Schoenberg correspondence for positive definite kernels) appeared in \cite{Sch38}. 
Schoenberg correspondences play an important role in noncommutative probability. Schoenberg correspondences have been proved by Schürmann for convolution semigroups of sesquilinear forms on $*$-coalgebras and for linear functionals on $*$-bialgebras \cite{Sch85}. The Schoenberg correspondence for $*$-bialgebras is also a corner stone of the theory of quantum Lévy processes as developed in \cite{Sch93}; it describes \emph{generating functionals} of Lévy processes on a $*$-bialgebra as linear functionals on the $*$-bialgebra that are hermitian and whose restriction to the kernel of the counit is positive. Lévy processes on $*$-bialgebras have, by definition, increments which are independent with respect to tensor independence. In contrast with classical probability, it is a remarkable feature of noncommutative probability that there are other independences which can be used to define Lévy processes. A particularly nice class of independences are the \emph{universal independences}; noncommutative random variables are called independent if their noncommutative joint distribution (modelled on the free product of algebras) coincides with the product distribution with respect to a \emph{universal product} in the sense of Ben Ghorbal \& Sch{\"u}rmann \cite{BGSc05}. Universal products are defined through axioms which reflect desirable properties of the corresponding independence. As explained in \cite{BGSc05}, a universal independence yields a convolution for linear functionals on dual semigroups. In \cite{ScVo14} Schürmann \& Voß proved Schoenberg correspondences for universal independences. Instead of concluding the result directly from the axioms of universal products, they use Fock space constructions which have been established for specific universal products; this is possible due to Muraki's classification theorem \cite{Mur02,Mur13} which states that a (positive) universal independence has to be one of the five well-studied examples: tensor, free, Boolean, monotone, and antimonotone independence. 

There are generalizations of universal products in which no classification is available, leading to other classes of Lévy processes. This raises the question whether Schoenberg correspondence also holds in these cases. In this paper we address the multivariate setting of independences associated with $m$-$d$-universal products in the sense of Manzel \& Sch{\"u}rmann \cite{MaSc17} (i.e., universal products for $d$-tuples of $m$-dimensional noncommutative distributions); the most prominent examples are probably Bożejko \& Speicher's $c$-freeness \cite{BoSp96} (corresponding to $m=1, d=2$) and Voiculescu's bifreeness \cite{Voi14} (corresponding to $m=2, d=1$), but many more such examples have been studied (e.g., (free-)free-Boolean independence by Liu in \cite{Liu19,Liu18p} or bimonotone independence of type II by the author in \cite{G23} and by Gu, Skoufranis \& Hasebe in \cite{GHS20}). In \cite{GHU23} the author, Hasebe \& Ulrich constructed deformations leading to continuous families of independences for $m=2,d=1$. All the beforementioned examples are  \emph{positive} in the sense that the product of states (more precisely: \emph{restricted states}, see Section \ref{sec:positivity}) is again a state. Gu \& Skoufranis' c-bifree independence ($m=2,d=2$) with its derivates bi-Boolean independence and (with Hasebe) bimonotone independence of type I (both $m=2,d=1$) \cite{GuSk17,GuSk19,GHS20}, as well as Lachs' $r$-$s$-independence ($m=d=1$) \cite{GLa15, Lac15} are not positive.

In an attempt to classify positive and symmetric $2$-$1$-universal products, called \emph{two-faced universal products}, the author and Var\v{s}o \cite{GVa24} present examples of independences for which positivity is yet open; we call them the \emph{exceptional} universal products. Knowing that Schoenberg correspondence holds also in this case is of particular interest. On the one hand it provides additional necessary criteria for positivity: certain Hankel determinants related to one parameter semigroups have to be positive. In his PhD thesis \cite{Var21}, Var\v{s}o calculated some of these Hankel determinants and could not find a contradiction to positivity, thus strengthening our belief that they are indeed positive universal products. In \cite{GHU23}, all positive two-faced universal products which can be ``represented'' on the tensor product or free product of Hilbert spaces are classified. In the light of these results, the exceptional universal products, should they turn out to positive, are bound to have more complicated operator models than all other known examples. A general Schoenberg correspondence would then establishes the existence of L{\'e}vy processes with respect to these independences and might even help to find associated operator models. There is quite some interest in ``generalized Gaussian processes'' which interpolate between the Brownian motions associated with $1$-$1$-independences (most notably Bo\.zejko and Speicher's $q$-Brownian motion \cite{BoSp92}). Every positive two-faced universal product gives rise to such an interpolation. 

In this paper we reinvestigate Sch{\"u}rmann \& Vo{\ss}' proof of the Schoenberg correspondence for positive universal products \cite[Section 3]{ScVo14}, observe that a variation of their approach allows us to furnish a proof which does not rely on Muraki's classification, and simultaneously generalize the proof to the multivariate situation, where the problem was still open (cf.\ \cite[Remark 4.5]{MaSc17}). On the way we present concise and applicable statements about the \emph{Lachs functor}, a versatile and useful tool to study of universal products developed by Ben Ghorbal \& Sch{\"u}rmann \cite{BGSc02,BGSc05}, Lachs \cite{Lac15}, and Manzel \& Sch{\"u}rmann \cite{MaSc17}. 

The paper is organized as follows. 
In Section \ref{sec:prelims} we define and review some properties of coalgebras and dual semigroups. A coalgebra $C$ has a comultiplication $\Delta\colon C\to C\otimes C$ which allows us to define a convolution for linear functionals on $C$ based on the tensor product $\varphi_1*\varphi_2:=(\varphi_1\otimes \varphi_2)\circ\Delta$. Dual semigroups are defined similar to coalgebras, but a dual semigroup $B$ is an algebra and has its comultiplication $\Lambda$ take values in the free product $B\sqcup B$. We will also need a multivariate version of dual semigroups called \emph{$m$-faced} dual semigroups to formulate the multivariate Schoenberg correspondences. In Section \ref{sec:mult-univ-prod} we define (multivariate) universal products, which are ways to form a product of linear functionals $\varphi_1,\varphi_2$ on algebras $A_1,A_2$ to obtain a linear functional $\varphi_1\odot\varphi_2$ on the free product $A_1\sqcup A_2$ of the domains. In particular, a universal product gives rise to a convolution of linear functionals $\varphi_1,\varphi_2\colon B\to\mathbb C$ on a dual semigroup $(B,\Lambda)$, defined as $\varphi_1\star\varphi_2:=(\varphi_1\odot\varphi_2)\circ\Lambda$. In Section \ref{sec:conv-expon-lachs} we discuss a key tool to understand universal products and convolutions, the \emph{Lachs functor} $\mathcal L$, which allows us to express (multivariate) dual semigroup convolutions $\varphi_1\star\varphi_2$ as usual coalgebraic convolutions on a bialgebra $\mathcal L(B)$. We formulate precise statements about $\mathcal L$ in Theorem \ref{th:Lachs-functor} in the very efficient language of monoidal categories and monoidal functors as laid out for example in \cite{AgMa10}; for readers unfamiliar with those notions, in Remark \ref{rem:Lachs-functor}, we try to provide enough intuition on what the Lachs functor does to follow the rest of this paper, also the statements of \ref{cor:Lachs-functor} and \ref{cor:exp}, which describe the relationship between coalgebraic and dual semigroup convolutions and convolution exponentials, can be understood without background knowledge in monoidal categories, although their proofs become quite lengthy (cf. \cite[Theorems 3.4 and 4.6]{BGSc05} for the concrete proofs in the univariate case).
Section \ref{sec:positivity} is devoted to the main result of this paper, the Schoenberg correspondence for multivariate universal products, which characterizes the linear functionals $\psi$ on a $*$-$m$-faced dual semigroup $B$ for which the convolution exponentials $\exp_\odot(t\psi)$ are (restrictions of) states on the unitization of $B$ for all $t>0$. In Section \ref{sec:gener-brown-moti} we indicate a possible application to the construction of Brownian motions and Fock spaces.

Before we dive into the details, let us vaguely outline what is the idea behind the proof of Schoenberg correspondence in \cite{ScVo14} and what is the difference in this paper. Consider a structure with an underlying topological vector space $X$ and
\begin{itemize}
\item a notion of \emph{positivity}, i.e., a closed subset $X_+\subset X$
\item two \emph{multiplications} $\cdot_1,\cdot_2\colon X\times X\to X$ which preserve positivity
\item two \emph{exponentials} $\exp_1,\exp_2\colon X\to X$.
\end{itemize}
Under good circumstances, one can express one exponential via the other as
\begin{equation}
\exp_2(x)= \lim_{n\to\infty} \exp_1({\textstyle\frac{1}{n}} x) \cdot_2 \ldots \cdot_2 \exp_1({\textstyle\frac{1}{n}}x)=\lim_{n\to\infty} \left(\exp_1({\textstyle\frac{1}{n}}x)\right)^{\cdot_2 n}.\label{eq:idea}
\end{equation}
Now suppose that there is a subset $G\subset X$ of \emph{generators} such that $\exp_1(t g)$ is positive for all $g\in G$ and $t>0$. Then one can conclude from Equation~\eqref{eq:idea} that $\exp_2(tg)$ is positive for all $g\in G$ and $t>0$. Universal products lead to a lot of possibilities to make this idea precise. In \cite{ScVo14} (inspired by similar results of Sch{\"u}rmann, Skeide \& Volkwardt for bialgebras \cite{SSV10}), different comultiplications $\Lambda_1,\Lambda_2\colon B\to B\sqcup B$ lead to different convolutions and convolution exponentials on the dual space $B'$ such that Equation~\eqref{eq:idea} holds. Therefore, Schoenberg correspondence can be transferred from easier to more complicated dual semigroups; the universal product is fixed beforehand and one has to establish one starting point for each universal product, which is where the dependence on classification theorems stems from. In this paper, the different exponentials correspond to different universal products and we still obtain Equation~\eqref{eq:idea}. This allows us in the end to reduce all Schoenberg correspondences in the context of universal products to Sch{\"u}rmann's Schoenberg correspondence for linear functionals on $*$-bialgebras.

\section{Coalgebras and dual semigroups}
\label{sec:prelims}

A coalgebra is a vector space $C$ with a coassociative comultiplication $\Delta\colon C\to C\otimes C$ (i.e., $(\Delta\otimes\mathrm{id})\Delta=(\mathrm{id}\otimes\Delta)\Delta$) and a counit $\delta\colon C\to\mathbb C$ (i.e., $(\delta\otimes\mathrm{id})\Delta=\mathrm{id}=(\mathrm{id}\otimes\delta)\Delta$); we identify $C\equiv C\otimes\mathbb C\equiv \mathbb C\otimes C$ and $(C\otimes C)\otimes C\equiv C\otimes (C\otimes C)$ in the natural way. In other words, coalgebras are comonoids in the monoidal category $(\mathrm{Vec},\otimes,\mathbb C)$. 
The dual $C'$ of a coalgebra $C$ is a unital algebra with respect to the convolution product
\(f*g=(f\otimes g)\circ \Delta\) and with unit $\delta$.

Two key properties make coalgebras much easier to work with than dual semigroups:
\begin{itemize}
\item The tensor product and, therefore, coalgebraic convolution are bilinear. (Universal products and convolutions on dual semigroups are not!)
\item Every finite dimensional subspace of a coalgebra is contained in a finite dimensional subcoalgebra (this is known as \emph{fundamental theorem for coalgebras}, for a proof see e.g.\ \cite[Theorem 1.4.7]{DNR01}).
\end{itemize}
These properties alone allow to establish a good functional calculus on the convolution algebra. For convenience of the reader, we present here the specific instances of this calculus which we will need in the sequel. 

\begin{lemma}[{cf.\ \cite[Section 4]{SSV10}}]\label{main-lemma}
  Let $(C,\Delta, \delta)$ be a coalgebra, $\psi\in C'$.
  \begin{enumerate}
  \item For every $c\in C$, the evaluation of the exponential series
    \[\exp_C(\psi)(c)=\sum_{n=0}^\infty \frac{\psi^{* n} (c)}{n!}\]
    is absolutely convergent in $\mathbb C$.
  \item Suppose that $R_n,n\in \mathbb N$, are linear functionals on $C$ such that for each $c\in C$ there is a constant $K_c\in \mathbb R_+$ with $|R_n(c)|\leq \frac{1}{n^2}K_c$ for all $n\in\mathbb N$. Then, for every $c\in C$,
    \[\exp_C(\psi)(c)=\lim_{n\to\infty} \left(\delta + \frac{\psi}{n}+R_n\right)^{\ast n}(c).\] 
  \end{enumerate}
\end{lemma}

\begin{proof}
  We sketch the arguments given in \cite{SSV10}.
  Fix $c\in C$ and let $\operatorname{coalg}(c)$ be smallest subcoalgebra of $C$ which contains $c$. Then $\operatorname{coalg}(c)$ is finite dimensional, so we can choose an arbitrary norm on $\operatorname{coalg}(c)$, and this induces (operator) norms on $\operatorname{coalg}(c)'$ and $L(\operatorname{coalg}(c))$; we denote all three norms simply $\|\cdot\|$. Then $L(\operatorname{coalg}(c))$ is a finite-dimensional unital normed algebra. With the convolution operator $T_\psi:=(\mathrm{id}\otimes \psi)\circ\Delta\restriction{\operatorname{coalg}(c)}\in L(\operatorname{coalg}(c))$, for every polynomial $p$ we have $p(\psi)=\delta(p(T_\psi))$. Therefore,
  \[ \sum_{n=k}^\ell\left|\frac{\psi^{* n} (c)}{n!}\right|= \sum_{n=k}^\ell\left|\frac{\delta(T_\psi^n(c))}{n!}\right|\leq \|\delta\| \|c\| \sum_{n=k}^\ell\frac{ \|T_\psi\|^n}{n!} \xrightarrow{k,\ell\to\infty} 0.\]
  The second claim is a specialization of \cite[Lemma 4.2]{SSV10}. The idea is to prove a corresponding result for Banach algebras, estimate the norm of $T_{R_n}=(\mathrm{id}\otimes R_n)\circ\Delta\restriction{\operatorname{coalg}(c)}$, and conclude 
  \[\left(\delta + {\textstyle\frac{\psi}{n}}+R_n\right)^{\ast n}(c) =  \delta\left(\mathrm{id} + {\textstyle\frac{T_\psi}{n}} + T_{R_n}\right)^n(c)\xrightarrow{n\to\infty}\delta \exp(T_\psi)(c)=\exp_C{\psi}(c). \]
\end{proof}

In this article, \emph{algebra} always means an associative, not necessarily unital algebra over the complex numbers $\mathbb C$. The category of algebras with algebra homomorphisms is denoted $\mathrm{Alg}$. The free product of algebras is denoted $\sqcup$, indicating that this is the coproduct in the category of algebras. This turns $\mathrm{Alg}$ into a monoidal category with the trivial algebra $\{0\}$ as unit object. 

An \emph{$m$-faced algebra} is an algebra $A$ together with subalgebras $A^1,\ldots, A^m$  (the \emph{faces} of $A$) which freely generate $A$, i.e., the canonical homomorphism $A^1\sqcup\ldots \sqcup A^m\to A$ is an isomorphism; we indicate this by writing $A=A^1\sqcup\ldots \sqcup A^m$. An \emph{$m$-faced algebra homomorphism} $j\colon A\to B$ is an algebra homomorphism between $m$-faced algebras $A,B$ with $j(A^k)\subset B^k$. The category of $m$-faced algebras with $m$-faced algebra homomorphisms is denoted $\mathrm{Alg}_m$. We consider the free product of $m$-faced algebras again an $m$-faced algebra with faces $(A\sqcup B)^k:=A^k\sqcup B^k$.  Note that the free product of $m$-faced algebras is the coproduct in the category $\mathrm{Alg}_m$. Therefore, $(\mathrm{Alg}_m,\sqcup,\{0\})$ is a monoidal category for every $m\in\mathbb N$.

It should be clear from the context whether an upper index labels a face or indicates a power, but to avoid confusion, we will sometimes explicitly write $A^{\oplus d}$ instead of $A^d$ for the direct sum $A\oplus\ldots\oplus A$ of $d$ copies when $A$ is an $m$-faced algebra.

\begin{remark}
  Two-faced algebras are related to Voiculescu's \emph{pairs of faces} \cite[Def. 2.1]{Voi14} in the following way. First, note that Voiculescu's unitality assumptions are appropriate to study (bi-)freeness, but the non-unital framework is more fitting for a general framework to study independences, so it is better to compare with Liu's non-unital version of the definition \cite[Def. 3.1]{Liu19} (of course there is also a unital version of multifaced algebras). A pair of faces in a non-commutative probability space $\mathcal (\mathcal A,\phi)$ is by definition a pair of algebra homomorphisms $j^1\colon B^1\to\mathcal A, j^2\colon B^2\to\mathcal A$, whose common codomain is the non-commutative probability space. This can equivalently be described by one homomorphism $j=j^1\sqcup j^2\colon B^1\sqcup B^2 \to\mathcal A$ whose domain is the two-faced algebra $B:= B^1\sqcup B^2$ with the faces $B^1,B^2$. Second, recall that in non-commutative probability the analogue of a random variable (i.e., a measurable map whose domain is a probability space) is an algebra homomorphism whose codomain is a non-commutative probability space. So, pairs of faces are simply pairs of non-commutative random variables, and two-faced algebras are the natural domains of the homomorphisms which model such pairs. Analogously, $m$-faced algebras are the natural domains of homomorphisms which model $m$-tuples of non-commutative random variables.
\end{remark}

\begin{definition}\label{def:dual-semigroup}
  An $m$-faced \emph{dual semigroup} is an $m$-faced algebra $B$ together with an $m$-faced algebra homomorphism $\Lambda\colon B\to B\sqcup B$ (the \emph{comultiplication}) which is \emph{coassociative} and fulfills the \emph{counit property}, i.e.,
  \[(\Lambda\sqcup \mathrm{id}) \Lambda=(\mathrm{id}\sqcup\Lambda)\Lambda\quad\text{and} \quad(\mathrm{id}\sqcup 0)\Lambda =\mathrm{id}=(0\sqcup\mathrm{id})\Lambda\] (the counit property makes the unique linear map $0\colon B\to \{0\}$ a \emph{counit} for $\Lambda$). In other words, an $m$-faced dual semigroup is a comonoid $(B,\Lambda,0)$ in the monoidal category $(\mathrm{Alg}_m,\sqcup,\{0\})$. A 1-faced dual semigroup is simply called \emph{dual semigroup}.
\end{definition}

In order to make it easily visible whether a comultiplication takes values in the tensor product (coalgebras) or the free product (dual semigroups), we will consistently write $\Delta$ in the former and $\Lambda$ in the latter situation.

\begin{example}
  The standard example of a dual semigroup is the tensor algebra $T(V)$ over a vector space with \emph{primitive} comultiplication $\Lambda(v)=v_1+v_2\in T(V)\sqcup T(V)$, where $v_1,v_2$ denote copies of $v$ interpreted as elements of the first and second factor $T(V)\supset V$ in the free product, respectively. If $V=V^1\oplus \cdots \oplus V^m$ is decomposed as a direct sum of subspaces, then $T(V)=T(V^1)\sqcup\cdots\sqcup T(V^m)$ inherits a decomposition as a free product of subalgebras. Furthermore, $\Lambda(T(V^k))\subset T(V^k)\sqcup T(V^k)$ shows that $T(V)$ can be considered an $m$-faced dual semigroup with faces $(T(V))^k=T(V^k)\subset T(V)$.
\end{example}

One can also define ($m$-faced) dual semigroups in a framework of unital algebras. We briefly sketch the relation of such \emph{unital dual semigroups} with the dual semigroups defined in Definition~\ref{def:dual-semigroup}. Details (for $m=1$, but $m>1$ works analogously) can be found in \cite[Section 3]{BGSc05}.

Consider the category $\mathrm{uAlg}$ of unital algebras with unital algebra homomorphisms. For $\mathrm{uAlg}$, the coproduct is given by the \emph{unital free product} $\sqcup_1$ (i.e., the free product with identification of units) and the initial object is $\mathbb C$. Therefore, $(\mathrm{uAlg},\sqcup_1,\mathbb C)$ is a monoidal category. 

Let $\unitization\colon \mathrm{Alg}\to\mathrm{uAlg}$ denote the \emph{unitization functor}; for an algebra $A$, $\unitization A$ is the unital algebra which contains $A$ as an ideal of codimension 1 (even if $A$ had a unit beforehand!) and for an algebra homomorphism $j\colon A\to B$, $\unitization j\colon \unitization A\to\unitization B$ is the unique unital extension of $j$. This is a \emph{monoidal functor}, which means, roughly speaking, that one can identify $\unitization (A\sqcup B)\equiv \unitization A \sqcup_1 \unitization B$.  A \emph{unital dual semigroup} is a unital algebra $\mathcal B$ with unital algebra homomorphisms $\Lambda\colon \mathcal B\to\mathcal B\sqcup_1\mathcal B$ and $\lambda\colon \mathcal B\to \mathbb C$ such that coassociativity ($(\Lambda\sqcup_1\mathrm{id})\Lambda=(\mathrm{id}\sqcup_1\Lambda)\Lambda$) and the counit property ($(\lambda\sqcup_1\mathrm{id})\Lambda=\mathrm{id}=(\mathrm{id}\sqcup_1\lambda)\Lambda$) are fulfilled; in short, a unital dual semigroup is a comonoid in $(\mathrm{uAlg},\sqcup_1,\mathbb C)$. There is a one-to-one correspondence between dual semigroups and unital dual semigroups (up to isomorphism) given by the following prescriptions:

\begin{itemize}
\item $(\mathcal B,\Lambda, \lambda)\mapsto (B,\Lambda\restriction B,0)$ with $B:=\operatorname{kern}\lambda$; one can check that $\Lambda(B)\subset B\sqcup B$ (identified with a subset of $\mathcal B\sqcup_1 \mathcal B$) and the rest follows easily.
\item $(A,\Lambda,0)\mapsto (\unitization A,\unitization \Lambda, \unitization 0)$. 
\end{itemize}

This establishes an equivalence of the categories of dual semigroups and unital dual semigroups.

\begin{remark}
  The unital variant is closer to Voiculescu's original definition of dual groups in \cite{Voi87} (in the context of pro-C*-algebras), while the nonunital variant is typically preferred in the context of universal products (cf.\ \cite{BGSc05,MaSc17}); the reason is that, in general, universal products do not respect identification of units, which complicates the formulas and statements considerably. Both versions are $H_0$-algebras in the sense of Zhang \cite{Zha91}.

    We also chose to define dual semigroups in the nonunital framework to harmonize better with the definition of universal products in Section \ref{sec:mult-univ-prod}. However, in certain contexts, e.g., in Section \ref{sec:positivity} of this paper where we study positivity of linear functionals on dual semigroups, one should think of the nonunital algebras as ideals in their unitizations in order to get the right intuition.
\end{remark}

\section{Multivariate universal products}
\label{sec:mult-univ-prod}

As opposed to the coalgebraic setting, there is no canonical candidate for a convolution product for functionals on a dual semigroup; linear functionals aren't morphisms in  $\mathrm{Alg}_m$, so we can not simply replace the tensor product in the convolution formula by the free product. In this section we recall the concept of (multivariate) \emph{universal products}; introduced in \cite{BGSc02} and generalized in \cite{MaSc17}, universal products are designed as ideal replacements for the tensor product in the definition of convolution.

\begin{definition}
  A $m$-faced $d$-valued linear functional, $m$-$d$-functional for short, is a pair $(A,\varphi)$, where $A$ is an $m$-faced algebra and $\varphi\colon A\to\mathbb C^d$ a linear map.
  By $\mathcal{F}_{m,d}$ we denote the category of $m$-$d$-functionals, i.e.,
  \begin{itemize}
  \item an object of $\mathrm{\mathcal F}_{m,d}$ is an $m$-$d$-functional 
  \item a morphism $j\colon(A_1, \varphi_1)\to(A_2, \varphi_2)$ is an $m$-faced algebra homomorphism $j\colon A_1\to A_2$ with $\varphi_1=\varphi_2\circ j$.
  \end{itemize}
  The obvious forgetful functor which associates with each $m$-$d$-functional $(A,\varphi)$ its domain $A$ is denoted $\operatorname{dom}\colon \mathcal{F}_{m,d}\to \mathrm{Alg}_m$.  Usually, we will write $\varphi$ instead of $(A,\varphi)$ and $\operatorname{dom}(\varphi)$ to refer to $A$. Note that $\operatorname{dom}(j)=j$ for every morphism.

  By $\unitization  \varphi$ we denote the unique unital extension $\unitization  \varphi\colon \unitization  A\to \mathbb C^d$.
  (If $\varphi$ is also an algebra homomorphism, there are possible conflicts of notation, but it should be quite clear from the context whether we deal with objects of $\mathcal F_{m,d}$ or with morphisms of some category.)
\end{definition}

Of course, a linear map $\varphi\colon A\to \mathbb C^d$ is uniquely determined by its components $\varphi_k\colon A\to\mathbb C$, $k\in\{1,\ldots, d\}$.  Accordingly, we denote the set of all $m$-$d$-functionals on an $m$-faced algebra $A$ by $A'^{d}$. We will often make use of the \emph{transpose} of $\varphi\in A'^{d}$ defined as $\varphi^{\mathrm{tr}}\colon A^{\oplus d}\to\mathbb C, \varphi^{\mathrm{tr}}(a_1,\ldots, a_d):= \sum_{k=1}^d \varphi_k(a_k)$. Clearly, $\mathrm{tr}\colon  A'^{d}\to (A^{\oplus d})'$ is a linear isomorphism and the inverse is denoted by the same symbol $\mathrm{tr}\colon  (A^{\oplus d})'\to A'^{d}$. 

\begin{remark}
  The category $\mathcal F_{m,d}$ is (trivially) isomorphic to the category $\mathrm{AlgP}_{d,m}$ of \cite{MaSc17}. The different choice of notation is mainly made because the reader should not be encouraged to think of all linear functionals on algebras as analogues of probability measures; while this interpretation is quite adequate for certain functionals which encode moments of noncommutative random variables, cumulant functionals are also represented by objects of $\mathcal F_{m,d}$ and play an equally important, but very different role.     
\end{remark}

Next, we present the categorical definition of universal products, the traditional definition can be found in Remark \ref{rem:up-trad} right below.

\begin{definition}
  An $m$-$d$-universal product is a bifunctor $\odot\colon \mathcal{F}_{m,d}\times \mathcal{F}_{m,d} \to \mathcal{F}_{m,d}$ such that
  \begin{enumerate}[label={(UP\arabic*)}]
  \item The bifunctor ${\odot}$ turns $\mathcal{F}_{m,d}$ into a monoidal category with unit object $0\colon \{0\}\to\mathbb C^d$.
  \item \label{UP2} The functor $\operatorname{dom}\colon (\mathcal F_{m,d},\odot,0)\to (\mathrm{Alg}_m,\sqcup,\{0\})$ is strictly monoidal, i.e., $\operatorname{dom}(\varphi_1\odot \varphi_2)=\operatorname{dom}(\varphi_1)\sqcup \operatorname{dom}(\varphi_2)$ and $j_1\odot j_2=\operatorname{dom}(j_1\odot j_2)=j_1\sqcup j_2$. 
  \end{enumerate}  
\end{definition}

\begin{remark}\label{rem:up-trad}
  Using the canonical identifications $A\sqcup \{0\}\equiv A \equiv\{0\}\sqcup A$ and $(A_1\sqcup A_2)\sqcup A_3\equiv A_1\sqcup (A_2\sqcup A_3)$, the previous definition can be paraphrased in the following more traditional way: An $m$-$d$-universal product is a binary product operation for $d$-valued linear functionals on $m$-faced algebras which associates with functionals $\varphi_1,\varphi_2$ on $m$-faced algebras $A_1,A_2$, respectively, a functional $\varphi_1\odot\varphi_2$ on $A_1\sqcup A_2$ such that
  \begin{itemize}
  \item $(\varphi_1\circ j_1)\odot(\varphi_2\circ j_2)=(\varphi_1\odot\varphi_2)\circ (j_1\sqcup j_2)$ (universality)
  \item $(\varphi_1\odot\varphi_2)\odot\varphi_3 = \varphi_1\odot(\varphi_2\odot\varphi_3)$ (associativity)
  \item $\varphi\odot 0=\varphi = 0\odot\varphi$ (unit condition) 
  \end{itemize}
  Finally, using universality and the fact that the embeddings $A_k\hookrightarrow A_1\sqcup A_2$ can be written $\mathrm{id}\sqcup 0 \colon A_1\sqcup \{0\}=A_1\to A_1\sqcup A_2$, $0\sqcup \mathrm{id}\colon \{0\}\sqcup A_2=A_2\to A_1\sqcup A_2$, the last condition is easily seen to be equivalent to
  \begin{itemize}
  \item $(\varphi_1\odot\varphi_2)\restriction A_i=\varphi_i$ (restriction property).
  \end{itemize} 
\end{remark}

\begin{remark}
  Working without condition \ref{UP2} is possible. However, as an effect of the universal property of the free product (being the coproduct in the category of $m$-faced algebras), this would not be a big generalization. For some details on this, we refer the reader to \cite[Remark 5.2]{GLS22}.
\end{remark}

\begin{definition}
  Let $(B,\Lambda)$ be an $m$-faced dual semigroup and $\odot$ an $m$-$d$-universal product. The \emph{convolution} of $d$-valued functionals $\varphi_1,\varphi_2\colon B\to\mathbb C^d$ is defined as
  \[\varphi_1\star \varphi_2:=(\varphi_1\odot\varphi_2)\circ \Lambda.\]
\end{definition}

\begin{remark}
  Note that the convolution depends on the choice of a universal product. Also be aware that universal products and, therefore, the just defined convolution are typically not bilinear. In particular, the $d$-valued functionals form a monoid, but not an algebra. It is also not meaningful to define a convolution exponential directly via the exponential series with respect to the operation $\star$. 
\end{remark}

Just as we use different symbols $\Delta,\Lambda$ for comultiplications into the tensor or free product, we use two different symbols $*,\star$ for the corresponding convolutions.

\section{Convolution exponentials and the Lachs functor}
\label{sec:conv-expon-lachs}

In this section we recall a key technique for dealing with convolution on dual semigroups developed by Ben Ghorbal \& Sch{\"u}rmann  \cite{BGSc02,BGSc05}, Lachs \cite{Lac15} and Manzel \& Sch{\"u}rmann \cite{MaSc17}. Remark \ref{rem:Lachs-functor} below explains their contributions and sheds some light on the underlying ideas. The results presented in this section, Theorem \ref{th:Lachs-functor}, Corollary \ref{cor:Lachs-functor}, and Corollary \ref{cor:exp}, were known before. However, the ideas are somewhat scattered in the references and not always formulated in the degree of generality we need, so we compress here the essence of \cite[Section 5.2]{Lac15} and \cite[Section 5]{MaSc17} in a way suitable for application in Section \ref{sec:positivity} as well as for possible future reference. 


Fix $m,d\in\mathbb N$. Let $\operatorname{uAlg}$ denote the category of unital algebras and $\mathrm{u}\mathcal F$ the category of linear functionals on unital algebras (the functionals need not be unital). Given an $m$-faced algebra $A$, we denote by $\mathcal L_{\mathrm{Alg}}(A)=S(A^{\oplus d})$ the symmetric tensor algebra over the direct sum of $d$ copies of $A$. Clearly, this defines a functor $\mathcal L_{\mathrm{Alg}}\colon \mathrm{Alg}_{m}\to \mathrm{uAlg}$ with the obvious action on morphisms $\mathcal L_{\mathrm{Alg}}(j)=S(j^{\oplus d})$. Note also that $\mathcal L_{\mathrm{Alg}}(\{0\})=\mathbb C$. Recall that the transpose of $\varphi\colon A\to\mathbb C^d$ is a map $\varphi^{\mathrm{tr}}\colon A^{\oplus d}\to \mathbb C$. We introduce a second functor, $\mathcal L_{\mathcal F}\colon \mathcal F_{m,d} \to \mathrm{u}\mathcal{F}$ with $\mathcal L_{\mathcal F}(\varphi)=S(\varphi^{\mathrm{tr}})$ and the same action on morphisms as $\mathcal L_{\mathrm{Alg}}$. 

Recall that a \emph{colax monoidal functor} (or \emph{cotensor functor}) is a functor $F\colon C_1\to C_2$ between monoidal categories $(C_k,\otimes_k,E_k)$ together with a natural transformation $F(\cdot\otimes_1 \cdot)\to F(\cdot)\otimes_2 F(\cdot)$ and a morphism $F(E_1)\to E_2 $ which fulfill certain compatibilities (\emph{coassociativity} and \emph{counitality}), see \cite[Definition 3.2]{AgMa10}.  

\begin{theorem}[cf.\ {\cite[Rem.~5.1]{MaSc17}}, {\cite[Thm.~5.2.4 + discussion]{Lac15}}]\label{th:Lachs-functor}
  Let $\odot$ be an $m$-$d$-universal product.
  Then, for all pairs of $m$-faced algebras $A_1,A_2$, there are unital algebra homomorphisms $\sigma^\odot_{A_1,A_2}\colon \mathcal L_{\mathrm{Alg}}(A_1\sqcup A_2)\to \mathcal L_{\mathrm{Alg}}(A_1)\otimes \mathcal L_{\mathrm{Alg}}(A_2)$ such that, with
  \[\sigma^\odot_{\mathrm{Alg}}:=(\sigma^\odot_{A_1,A_2})_{A_1,A_2\in \mathrm{Alg}_{m}},\qquad \sigma^\odot_\mathcal F:=(\sigma^\odot_{\operatorname{dom}\varphi_1,\operatorname{dom}\varphi_2})_{\varphi_1,\varphi_2\in\mathcal F_{m,d}}\]
  the following holds:
  \begin{itemize}
  \item $(\mathcal L_{\mathrm{Alg}},\sigma^\odot_{\mathrm{Alg}},\mathrm{id}_{\mathbb C})\colon(\mathrm{Alg}_{m},\sqcup,\{0\})\to(\mathrm{uAlg},\otimes, \mathbb C)$ is a colax monoidal functor.
  \item $(\mathcal L_{\mathcal F},\sigma^\odot_{\mathcal F},\mathrm{id}_{\mathbb C})\colon(\mathcal {F}_{m,d},\odot,\{0\})\to(\mathrm{u}\mathcal F,\otimes, \mathbb C)$ is a colax monoidal functor, in particular, for all pairs of $m$-$d$ functionals $\varphi_1$, $\varphi_2$, it holds that 
    \begin{align}
      S\bigl((\varphi_1\odot\varphi_2)^{\mathrm{tr}}\bigr) = \bigl(S(\varphi_1^{\mathrm{tr}}) \otimes S(\varphi_2^{\mathrm{tr}})\bigr)\circ \sigma^\odot_{A_1,A_2}.\label{eq:Lachs-functor}
    \end{align}
  \end{itemize}
\end{theorem}

\begin{definition}
  The \emph{Lachs functors} are the two colax monoidal functors
  \[\mathcal L^{\odot}_{\mathrm{Alg}}:=(\mathcal L_{\mathrm{Alg}},\sigma^\odot_{\mathrm{Alg}},\mathrm{id}_{\mathbb C}),\qquad\mathcal L_{\mathcal F}^{\odot}:=(\mathcal L_{\mathcal F},\sigma^\odot_{\mathcal F},\mathrm{id}_{\mathbb C}).\]
  (The natural transformations $\sigma^\odot_{\mathrm{Alg}}$ and $\sigma^\odot_{\mathcal F}$, which depend on $\odot$, are an integral part of the structure of the Lachs functors.)
\end{definition}

\begin{remark}
  Note that $\operatorname{dom}(\mathop{\mathcal L_{\mathcal F}}\varphi)=\mathcal L_{\mathrm{Alg}}(\operatorname{dom}\varphi)$ (where on the left hand side $\operatorname{dom}$ is the forgetful functor $\mathrm{u}\mathcal{F}\to \mathrm{uAlg}$ that maps each unital functional to its domain). Furthermore, the algebra homomorphisms which constitute  $\sigma^\odot_{\mathcal F}$ only depend on the domains of the functionals and then agree with the corresponding components of $\sigma^\odot$. So $\mathcal L^{\odot}_{\mathrm{Alg}}$ is, in this sense, the colax monoidal functor $(\mathrm{Alg}_{m},\sqcup,\{0\})\to(\mathrm{uAlg},\otimes, \mathbb C)$ induced by $\mathcal L^\odot_{\mathcal F}$ under taking domains. On the other hand,  $\mathcal L^{\odot}_{\mathcal{F}}$ is the canonical lift of $\mathcal L^{\odot}_{\mathrm{Alg}}$ to a colax monoidal functor $(\mathcal {F}_{m,d},\odot,\{0\})\to(\mathrm{u}\mathcal F,\otimes, \mathbb C)$.  This somehow justifies the singular expression ``the Lachs functor'' to refer to either or both of these two functors.
  In the sequel, we will sometimes drop the subscripts and simply write $\mathcal L^\odot$ and $\mathcal L$.
\end{remark}

Recall that a bialgebra is a vector space with compatible structures of coalgebra and unital algebra (i.e., the comultiplication and counit are unital algebra homomorphisms). In the categorical language, a bialgebra is a comonoid in $(\mathrm{uAlg},\otimes, \mathbb C)$. 

\begin{corollary}[cf.\ {\cite[Thm.~5.1 + Rem.~5.1]{MaSc17}}, {\cite[Prop.~5.2.1 + Eq.~(5.8)]{Lac15}}]\label{cor:Lachs-functor}
  The Lachs functor $\mathcal L^\odot_{\mathrm{Alg}}$ induces a functor between the categories of comonoids, i.e., if $B$ is an $m$-faced dual semigroup with comultiplication $\Lambda$, then the unital algebra $\mathcal L(B)=S(B^{\oplus d})$ is a bialgebra, with respect to the comultiplication $\sigma^\odot_{B,B}\circ S(\Lambda^d)$ and counit $S(0)$; this bialgebra is denoted $\mathcal L^\odot(B)$, siginifying that the colax monoidal structure of the Lachs functor is used and the resulting bialgebra depends on the universal product $\odot$. Furthermore, $\mathcal L_{\mathcal F}$ intertwines the convolutions $\star$ on $B'$ and $*$ on $\mathcal L^\odot(B)'$, i.e., 
  \[\mathcal L(\varphi_1\star\varphi_2)=S(\varphi_1\star\varphi_2)^{\mathrm{tr}}=S(\varphi_1^{\mathrm{tr}})*S(\varphi_2^{\mathrm{tr}})=\mathcal L(\varphi_1)*\mathcal L(\varphi_2).\]
\end{corollary}

\begin{proof}
  Colax monoidal functors preserve comonoids \cite[Proposition 3.29]{AgMa10}. That $\mathcal L$ intertwines convolutions is a direct consequence of Equation \eqref{eq:Lachs-functor}.
\end{proof}

\begin{remark}\label{rem:Lachs-functor}
  The construction underlying the Lachs functor dates back to Ben Ghorbal \& Schürmann \cite{BGSc02,BGSc05}. Its functorial properties were exhibited by Lachs \cite[Section 5.2]{Lac15}. Manzel \& Sch{\"u}rmann generalized the construction to $m$-$d$-universal products in \cite[Section 5]{MaSc17}. Theorem~\ref{th:Lachs-functor} is closely related to the ``central structural theorem'' for universal products, which states, roughly speaking, that an $m$-$d$-universal product $\varphi_1\odot\varphi_2(a_1\ldots a_n)$ can always be written as a polynomial in ``sub-evaluations'' $\varphi_i(a_{k_1}\ldots a_{k_{\ell}})$ with coefficients which are independent of the $\varphi_i$ and $a_k$ (cf.\ \cite[Theorem 4.2]{MaSc17} for the precise statement). We are satisfied with showing one example which hopefully makes this relation apparent for the reader. Let $\odot_f$ be the free product of linear functionals used in free probability. Then
\begin{align*}
\varphi_1\odot_f \varphi_2(a_1b_1a_2b_2)
  &=\quad\varphi_1(a_1a_2)\varphi_2(b_1)\varphi_2(b_2)\\
  &\quad+\varphi_1(a_1)\varphi_1(a_2)\varphi_2(b_1b_2)\\
  &\quad-\varphi_1(a_1)\varphi_1(a_2)\varphi_2(b_1)\varphi_2(b_2) 
\end{align*}
for all $\varphi_1,\varphi_2$ and all $a_1,a_2\in\operatorname{dom}\varphi_1, b_1,b_2\in\operatorname{dom}\varphi_2$. Therefore, a linear map $\sigma$ with
\begin{align*}
\sigma(a_1b_1a_2b_2)
  &=\quad a_1a_2 \otimes (b_1\cdot b_2) \\
  &\quad+(a_1\cdot a_2)\otimes b_1b_2\\
  &\quad-(a_1\cdot a_2)\otimes (b_1\cdot b_2) 
\end{align*}
(``$\cdot$'' the multiplication of the symmetric tensor algebra) fulfills Equation~\eqref{eq:Lachs-functor} when evaluated at $a_1b_1a_2b_2\in \operatorname{dom}\varphi_1 \sqcup \operatorname{dom}\varphi_2$. By the central structural theorem, this scheme works for all functionals $\varphi_1,\varphi_2$ and all elements of $\operatorname{dom}\varphi_1 \sqcup \operatorname{dom}\varphi_2$ and leads to Theorem \ref{th:Lachs-functor}. One can also go the other way round, i.e., the natural transformation of the Lachs functor determines a ``polynomial form'' of the corresponding universal product via Equation~\eqref{eq:Lachs-functor} and thus implies the central structural theorem.
\end{remark}

Consider a bialgebra $B$ with multiplication $\mu$, comultiplication $\Delta$ and counit $\delta$. The exponential $\exp_*$ on the dual $B'$ maps a \emph{derivation} $d$ (with respect to the  counit, i.e.\ $d(ab)=d(a)\delta(b)+\delta(a)d(b)$) to a \emph{character} $\exp_*(d)$ (a multiplicative linear functional); indeed,
\begin{align*}
  d^{* n}\circ \mu=d^{\otimes n}\circ \Delta^{(n)}\circ \mu=d^{\otimes n}\circ \mu^{\otimes n}\circ \Delta_{B\otimes B}^{(n)}=(d\circ\mu)^{*n}=(d\otimes\delta + \delta\otimes d)^{*n}
\end{align*}
and therefore, since $d\otimes\delta$ and $\delta\otimes d$ commute under the convolution on $(B\otimes B)'$,
\begin{multline*}
  \exp_* d\circ \mu = \sum \frac{1}{n!} d^{*n} \circ\mu = \exp_*(d\otimes\delta + \delta\otimes d)=(\exp_* d\otimes\delta)* (\delta\otimes \exp_* d)\\=\exp_* d\otimes\exp_* d.
\end{multline*}

For a linear functional $f$, $S(f)$ is the unique extension of $f$ to a character on $S(\operatorname{dom}f)$. We denote by $D(f)$ the unique extension of $f$ to an $S(0)$-derivation on $S(\operatorname{dom}f)$, i.e.,
\[D(f)(a_1\cdot\ldots \cdot a_n)=
  \begin{cases}
    f(a_1) & n=1,\\
    0 & n\neq 1.
  \end{cases}
\]
This gives us yet another functor, the \emph{differential Lachs functor}, $\mathcal D\colon \mathcal F_{m,d}\to \mathrm{u}\mathcal F$ with $\mathcal D(\psi)=D(\psi^{\mathrm{tr}})$. With this notation, we can comfortably describe convolution exponentials on $m$-faced dual semigroups. 

\begin{corollary}[Generalizes corresponding statements in {\cite[Rem.~7.1]{MaSc17}} and {\cite[Prop.~5.2.5 + discussion]{Lac15}}]\label{cor:exp}
  For an $m$-$d$-universal product $\odot$ and an $m$-faced dual semigroup $B$ there exists a unique map $\exp_\odot \colon B'^d\to B'^d$ such that, for all $\psi\in \mathcal F_{m,d},\operatorname{dom}\psi=B$,
  \[\mathcal L(\exp_{\odot}(\psi))=\exp_{\mathcal L^\odot(B)} \mathcal D(\psi)\]
  or, equivalently,
  \[\exp_\odot(\psi)=(\exp_{\mathcal L^\odot(B)} \mathcal D(\psi)\restriction B^{\oplus d})^{\mathrm{tr}}. \]
\end{corollary}

\begin{proof}
  Taking $\exp_\odot(\psi):=(\exp_{\mathcal L^\odot(B)} \mathcal D(\psi)\restriction B^{\oplus d})^{\mathrm{tr}}$ as a definition, and writing $\exp_{*}:=\exp_{\mathcal L^\odot(B)}$ for the coalgebraic convolution exponential for short, we see that
  \[
    \mathcal L(\exp_\odot \psi)
    =\mathcal L((\exp_{*} \mathcal D(\psi)\restriction B^{d})^{\mathrm{tr}})
    =S(\exp_{*} \mathcal D(\psi)\restriction B^{d})
    =\exp_{*} \mathcal D(\psi),
  \]
  where the last equality follows from the fact that $\exp_{*} \mathcal D(\psi)$ is a character and, trivially, an extension of its restriction $\exp_{*} \mathcal D(\psi)\restriction B^{\oplus d}$.
  Uniqueness follows from restricting to $B^{\oplus d}$.
\end{proof}

\begin{example}\label{ex:Lachs-functor}
  To digest what is going on, let us look at a well-studied example, the ``Standard Gaussian'' cumulant functional and its exponential with respect to tensor and Boolean product (without giving all the details, which can be found for example in \cite[Section 7.2]{Lac15}). In this example, we only deal with 1-1-universal products, i.e.~$m=d=1$. Let $B:= T(\mathbb C x)$, the tensor algebra over a one-dimensional vector space $\mathbb C x$ (which is just the algebra of polynomials in one indeterminate $x$ which vanish at $0$). The set $\langle x\rangle:=\{x^n:n\in\mathbb N\}$ is a basis of $B$. Therefore, $\mathcal L(B)$ is the polynomial algebra $\mathbb C[\langle x\rangle]=\mathbb C[x^1,x^2,\ldots]$ in infinitely many commuting(!) indeterminates $x^n, n\in\mathbb N$ (the internal multiplication of $B$ is forgotten). Recall that we denote by $x_1,x_2$ the element $x$ in the first and second copy of $B$ inside $B\sqcup B$, respectively, and note that as a vector space $B\sqcup B$ is spanned by elements of the form $x_1^{t_1}x_2^{s_1}\ldots x_1^{t_k} x_2^{s_k}$, where the exponents are natural numbers, all except the first and last exponent non-zero, and at least one exponent non-zero (an exponent of $0$ means omission of the factor because $B$ is not assumed unital).  We consider $B$ a dual semigroup with the comultiplication $\Lambda\colon B\to B\sqcup B$ determined by $\Lambda(x)=x_1+x_2$. Therefore, for each a 1-1-universal product $\odot$, the algebra $\mathcal L(B)$ is turned into a bialgebra $L^\odot(B)$.  First assume that $\odot=\otimes$ is the tensor product. Then $\sigma^\otimes_{B,B}\colon S(B\sqcup B)\to S(B)\otimes S(B)$ is the unique unital algebra homomorphism which is given on elements of $B\sqcup B$ by
  \[\sigma^\otimes_{B,B}(x_1^{t_1}x_2^{s_1}\ldots x_1^{t_k} x_2^{s_k})=x^{t_1+\ldots +t_k}\otimes x^{s_1+\ldots+s_k},\]
  in accordance with
  \[\varphi_1\otimes\varphi_2(x_1^{t_1}x_2^{s_1}\ldots x_1^{t_k} x_2^{s_k})=\varphi_1(x^{t_1+\ldots +t_k}) \varphi_2(x^{s_1+\ldots+s_k}).\]
  (In this case, $\sigma^\otimes_{B,B}(B\sqcup B)\subset B\otimes B$.)
  $\mathcal L(B)$ becomes the bialgebra $\mathcal L^\otimes(B)$ with comultiplication
  \begin{align*}
    \Delta^\otimes (x^n)=\sigma^\otimes_{B,B}(x_1+x_2)^n=\sum_{k=0}^n \binom{n}{k} x^k\otimes x^{n-k}.
  \end{align*}
  For the \emph{Gaussian cumulant functional}, 
  \[\psi(x^{n}):=
    \begin{cases}
      1 & \text{if $n=2$}\\
      0 & \text{otherwise}
    \end{cases}
  \]
  one obtains
  \begin{align*}
    \MoveEqLeft\exp_{\mathcal L^{\otimes}(B)}(D(\psi))(x^n)=\sum_{r=0}^{\infty} \frac{1}{r!}D(\psi)^{* r}(x^n)\\&=\sum_{r=0}^{\infty} \frac{1}{r!} \binom{n}{k_1,\ldots,k_r}D(\psi)(x^{k_1})\ldots D(\psi)(x^{k_r})\\
    &=
      \begin{cases}
        \frac{1}{(n/2)!}\binom{n}{2,\ldots,2}=(n-1)!! &\text{if $n$ is even,}\\
        0&\text{otherwise}
      \end{cases}
  \end{align*}
  the moments of the standard Gaussian distribution.

  For $\odot = \diamond$ the Boolean product, one has
  \[\sigma_{B,B}^\diamond(x_1^{t_1}x_2^{s_1}\ldots x_1^{t_k} x_2^{s_k})= (x^{t_1}\cdot\ldots\cdot x^{t_k})\otimes (x^{s_1}\cdot\ldots\cdot x^{s_k}) \]
  in accordance with
  \[\varphi_1\diamond\varphi_2(x^{t_1+\ldots +t_k}\otimes x^{s_1+\ldots+s_k})= \varphi_1(x^{t_1})\ldots \varphi_1(x^{t_k}) \varphi_2(x^{s_1})\ldots \varphi_2(x^{s_k})\]
  (of course, $\varphi_1(x^{t_1})$ and $\varphi_2(x^{s_k})$ are to be omitted if $t_1=0$ or $s_k=0$, respectively).
  The comultiplication of $\mathcal L^{\diamond}(B)$ is given by
  \[\Delta^{\diamond}(x^n)=\sigma^{\diamond}_{B,B}(x_1+x_2)^n\]
  For $n=4$ this results in
  \begin{align*}
    \Delta^{\diamond}(x^4)
    &=\begin{multlined}[t][.8\textwidth]
      \sigma^\diamond_{B,B}(x_1^4 + x_1^3x_2 + x_2x_1^3 +x_1^2x_2x_1+x_1x_2x_1^2 \\
      + x_1^2x_2^2+x_2^2x_1^2+x_1x_2^2x_1 + x_2x_1^2x_2 + x_1x_2x_1x_2
      +x_2x_1x_2x_1 \\
      + x_1x_2^3 + x_2^3 x_1 + x_2x_1x_2^2 + x_2^2x_1x_2 + x_2^4)
                            \end{multlined}\\
    &=
      \begin{multlined}[t][.8\linewidth]
        x^4\otimes 1 + 2\, x^3\otimes x + 2\,(x\cdot x^2)\otimes x \\
        + 2\,x^2\otimes x^2 + (x\cdot x)\otimes x^2 + x^2\otimes(x\cdot x) + 2\,(x\cdot x)\otimes (x\cdot x) \\
        + 2\, x\otimes x^3 + 2\,x\otimes (x\cdot x^2)  + 1\otimes x^4.
      \end{multlined}
  \end{align*}
  For the same Gaussian cumulant functional $\psi$ as before, one obtains
    \begin{align*}
      \MoveEqLeft\exp_{\mathcal L^{\diamond}(B)}(D(\psi))(x^4)=\sum_{r=0}^{\infty} \frac{1}{r!}D(\psi)^{* r}(x^4)\\
      &=\frac{1}{2!} (D(\psi)\otimes D(\psi))\Delta^{\diamond}(x^{4})=1
    \end{align*}
    With careful analysis of the expressions, one sees that in fact for all $n\in\mathbb N$,
     \begin{align*}
       \MoveEqLeft\exp_{\mathcal L^{\diamond}(B)}(D(\psi))(x^n)=\sum_{r=0}^{\infty} \frac{1}{r!}D(\psi)^{* r}(x^n)\\
       &=
      \begin{cases}
        \frac{(n/2)!}{(n/2)!} D(\psi)^{\otimes n/2} (x^2\otimes\ldots\otimes x^2)=1 &\text{if $n$ is even,}\\
        0&\text{otherwise}
      \end{cases}
  \end{align*}
  are the moments of Bernoulli distribution $\frac{1}{2}\delta_{-1} + \frac{1}{2}\delta_1$ which appears in the central limit theorem for Boolean independence.

  Note that in both cases $\exp_\odot(\psi)=\exp_{\mathcal{L^\odot}(B)}(D(\psi))|_B\neq\psi$ and, therefore, the character $\exp_{\mathcal{L^\odot}(B)}(D(\psi))$ is not just $S(\psi)$. 
  
\end{example}
\section{Positivity}
\label{sec:positivity}

For applications in noncommutative probability, it is desirable that universal products can be formed for distributions of noncommutative random variables. There is a small caveat here because positivity is best expressed in a unital framework -- many examples of universal products are, however, not fully compatible with the unital structures. The way out is to only consider unital $*$-algebras $\mathcal A$ of the form $\mathcal A=\unitization A$ for a given $*$-algebra $A$. (Alternatively, what is practically the same, one assumes that $\mathcal A$ is an \emph{augmented $*$-algebra}, i.e., a unital algebra with a $*$-character $\varepsilon\colon \mathcal A\to \mathbb C$; then, with $A:=\operatorname{kern}\varepsilon$, there is a canonical isomorphism $\mathcal A\cong \unitization A$.)

In the following definition, we write $x\geq0$ for a vector $x\in\mathbb C^d$ if all its components are non-negative real numbers.
\begin{definition}
A $d$-valued linear functional $\varphi\colon A\to\mathbb C^d$ on a $*$-algebra $A$ is called
\begin{itemize}
\item a \emph{restricted state} if its unitization $\unitization  \varphi$ is positive on the unital $*$-algebra $\unitization  A$, i.e., if
  \[\unitization \varphi(a^*a)\geq0 \text{ for all $a\in \unitization  A$},\]
\item a \emph{restricted generating functional} if it is hermitian and fulfills
  \[\varphi(a^*a)\geq0 \text{ for all $a\in A$}.\]
\end{itemize}
\end{definition}

For brevity, we decided to omit ``$d$-valued'' when it comes to (restricted) states and generating functionals, but note that $\varphi$ is a restricted state if and only if all \emph{components} of $\unitization \varphi$ are states.

\pagebreak[3]
\begin{remark}\leavevmode\makeatletter\@beginparpenalty\@M\makeatother
  \begin{itemize}
  \item Restricted states are automatically hermitian. Indeed, a restricted state $\varphi\colon A\to\mathbb C$ induces a positive sesquilinear form on $\unitization A$, $(a,b)_\varphi:=\varphi(a^{*}b)$, therefore $\varphi(a^{*})=(a,1)_\varphi=\overline{(1,a)_\varphi}=\overline{\varphi(a)}$. Note that positivity on $A$ is not sufficient for the argument to work.
  \item A restricted generating functional $\psi$ extends to a \emph{generating functional} $\hat\psi\colon\unitization \operatorname{dom}\psi \to\mathbb C^d$, i.e., to a $\hat\psi$ which is hermitian, vanishes at 1, and is positive on the ideal $\operatorname{dom}\psi\subset \unitization \operatorname{dom}\psi$ (\emph{conditionally positive}).
  \end{itemize}
\end{remark}

The free product of $*$-algebras carries a natural involution (the unique one extending the involutions of the factors) and will thus be considered a $*$-algebra.

\begin{definition}
An $m$-faced $*$-algebra is an $m$-faced algebra with an involution such that the faces are $*$-subalgebras. A $*$-$m$-faced dual semigroup is an $m$-faced dual semigroup with an involution such that the faces are $*$-subalgebras and the comultiplication is a $*$-homomorphism (i.e, a comonoid in the category of $m$-faced $*$-algebras). 
\end{definition}
\begin{definition}
  An $m$-$d$-$*$-distribution is a pair $(A,\varphi)$, where $A$ is an $m$-faced $*$-algebra and $\varphi\colon A\to \mathbb C^d$ is a restricted state.
\end{definition}

By simply ignoring the $*$-structure, we can lift a universal product to the category $*\text{-}\mathcal F_{m,d}$ of $d$-valued linear functionals on $m$-faced $*$-algebras, leaving us with the question how much of the inherent positivity structure is respected by a given universal product.

\begin{definition}
  Let $\odot$ be an $m$-$d$-universal product. We say that $\odot$ is
  \begin{itemize}
  \item \emph{positive} if the universal product $\varphi_1\odot\varphi_2$ of $m$-$d$-distributions $\varphi_1,\varphi_2$ is an $m$-$d$-distribution,
  \item \emph{half positive} if the universal product powers $\varphi^{\odot n}$ of an $m$-$d$-distribution $\varphi$ are $m$-$d$-distributions for all $n\in\mathbb N$,
  \item \emph{Schoenberg} if the convolution exponential $\exp_\odot(\psi)$ is a restricted state for every restricted generating functional $\psi$ on an $*$-$m$-faced dual semigroup.
  \end {itemize}
\end{definition}

If $\psi$ is a restricted generating functional on a $*$-algebra, then also $t\psi$ is a restricted generating functional for all $t>0$; therefore, if $\odot$ is Schoenberg and $\psi$ is a restricted generating functional on a $*$-$m$-faced dual semigroup $B$, then $\exp_\odot t\psi$ are restricted states for all $t>0$. A standard derivation argument can be used to show that if $\exp_\odot t\psi$ is a restricted state for all $t>0$, then $\psi$ must be a restricted generating functional: indeed, for $b\in B$, the map $f\colon \mathbb R_+\to \mathbb C, t\mapsto \exp_\odot t\psi(b^*b)$ is differentiable, positive and vanishes at zero, therefore $\psi(b^*b)=\left.\frac{\mathrm d}{\mathrm{d}t} f(t)\right|_{t=0}\geq0$. Combining the last two arguments, we obtain the following characterization of the Schoenberg property. 
\begin{observation}
  A universal product is Schoenberg if and only if the \emph{Schoenberg correspondence} holds:
  \[\psi \text{ restricted generating functional} \iff \exp_\odot{t\psi} \text{ restricted state for all $t>0$}\]
\end{observation}

\begin{observation}
  Let $B$ be a $*$-$m$-faced dual semigroup and $\odot$ an $m$-$d$-universal product. Since the comultiplication $\Lambda$ is a $*$-homomorphism, we have the following implications:
  \begin{itemize}
  \item If $\odot$ is positive, then $\varphi_1\star\varphi_2=(\varphi_1\odot\varphi_2)\circ\Lambda$ is a restricted state for all restricted states $\varphi_1,\varphi_2\colon B\to\mathbb C^d$.
  \item If $\odot$ is half positive, then $\varphi^{\star n}=\varphi^{\odot n}\circ \Lambda_n$ is a restricted state for every restricted state $\varphi\colon B\to\mathbb C^d$, where $\Lambda_n\colon B\to B^{\sqcup n}$ denotes the iterated comultiplication.
\end{itemize}

\end{observation}

One might think that the Schoenberg correspondence for dual semigroups follows trivially from Sch{\"u}rmann's Schoenberg correspondence for $*$-bialgebras \cite{Sch85} by application of the Lachs functor, but it is a bit more complicated than that: the natural transformation $\sigma$ of the Lachs functor (Theorem~\ref{th:Lachs-functor}) does not fulfill any obvious compatibility with the internal multiplication (in particular, $\sigma\restriction B^d\sqcup B^d\colon B^d\sqcup B^d\to S(B^d)\otimes S(B^d)$ is not a $*$-homomorphism), so there is no direct transfer of positivity properties! Still, it is \emph{a} key tool in the proof of Sch{\"u}rmann and Vo{\ss}' Schoenberg correspondence for positive universal products with $m=d=1$ and it is \emph{the} key tool in the proof we will present for the general $m$-$d$-case.

\begin{remark}
  Roughly speaking, the strategy of Sch{\"u}rmann \& Vo{\ss} \cite{ScVo14} is to apply Lemma~\ref{main-lemma} (2) in order to reduce the Schoenberg correspondence on a complicated dual semigroup to the Schoenberg correspondence on a dual semigroup with primitive comultiplication; using Muraki's classification, it is known that Schoenberg correspondence holds in those cases by explicit construction of a L{\'e}vy process on a suitable Fock space via quantum stochastic calculus. Our strategy is to apply the same Lemma~\ref{main-lemma} (2) in order to reduce the Schoenberg correspondence for a complicated universal product to the Schoenberg correspondence for the tensor product. This has the huge advantage that the proof is then independent of classification results and specialized quantum stochastic calculi, which are not available in the multivariate setting.
\end{remark}

\begin{theorem}\label{main-theorem}
  For an $m$-$d$-universal product $\odot$, the following implications hold:
  \begin{center}
    positive $\implies$ half positive $\implies$ Schoenberg
  \end{center}
\end{theorem}

\begin{proof}
  The first implication is trivial. For the second implication, let $\odot$ be an arbitrary half positive $m$-$d$-universal product and $\otimes$ the $m$-$d$-universal product given by componentwise tensor product (ignoring the faces), i.e.,
  \[(\varphi_1\otimes\varphi_2)_k=(\varphi_1)_k\otimes (\varphi_2)_k.\]
  Let $B$ be a $*$-$m$-faced dual semigroup and $\psi$ a restricted generating functional. We note that, for $\mathbf{b}\in B^{\oplus d}$,
  \begin{align*}
    \left|\exp_{\mathcal L^{\otimes}(B)}\left( \frac{1}{n} \mathcal D(\psi)\right)(\mathbf{b})- \delta(\mathbf{b})-\frac{1}{n}\mathcal D(\psi)(\mathbf{b})\right|
    \displaystyle \leq \frac{1}{n^2} \sum_{k=2}^\infty \frac{\bigl|\mathcal D(\psi)^{* k}(\mathbf{b})\bigr|}{k!}\leq \frac{1}{n^2}K_\mathbf{b}
  \end{align*}
  by Lemma~\ref{main-lemma} (1). Therefore, Lemma~\ref{main-lemma} (2) combined with Corollary~\ref{cor:exp}, for $\odot$ as well as for $\otimes$, implies
  \begin{align*}
    \mathcal L(\exp_\odot(\psi))&(\mathbf{b})
    =\exp_{\mathcal L^\odot(B)}(\mathcal D(\psi))(\mathbf{b})\\
    &=\lim_{n\to\infty} \left(\exp_{\mathcal L^{\otimes}(B)} \left(\textstyle\frac{1}{n}\mathcal D(\psi)\right)*_{\mathcal L^\odot(B)}\cdots *_{\mathcal L^\odot(B)}\exp_{\mathcal L^{\otimes}(B)}\left(\textstyle\frac{1}{n}\mathcal D(\psi)\right)\right)(\mathbf{b})\\
    &=\lim_{n\to\infty} \Bigl(\mathcal L\left(\exp_{\otimes}(\textstyle\frac{1}{n}\psi)\right)*_{\mathcal L^\odot(B)}\cdots *_{\mathcal L^\odot(B)}\mathcal L\left(\exp_{\otimes}(\textstyle\frac{1}{n}\psi)\right)\Bigr)(\mathbf{b})\\
    &=\lim_{n\to\infty} \mathcal L\left(\exp_\otimes(\textstyle\frac{1}{n}\psi)^{\star_\odot n}\right)(\mathbf{b}).
  \end{align*}
  Therefore, for $b\in B$, we find
  \[\exp_\odot(\psi)(b)=\lim_{n\to\infty} \left(\exp_\otimes(\textstyle\frac{1}{n}\psi)\right)^{\star_\odot n}(b).\]
  
  By half positivity, $\exp_\odot(\psi)$ is a restricted state if $\exp_\otimes(\textstyle\frac{1}{n}\psi)$ is a restricted state for all $n$. So we are done if we can show that $\otimes$ is Schoenberg. As the tensor product does not depend on the faces and does not mix the components, it is enough to show that it is Schoenberg for $m=d=1$. This follows for example from the results of \cite{ScVo14}. We present an alternative path. Note that the canonical $*$-homomorphism $\hat\sigma\colon B\sqcup B\to \unitization B\otimes \unitization B$ with $b_1\mapsto b\otimes 1, b_2\mapsto 1\otimes b$ coincides, as a linear map, with $\sigma_\otimes\restriction B\sqcup B$ from Theorem~\ref{th:Lachs-functor}. It follows that $\unitization  B\subset \mathcal S(B)$ is a $*$-bialgebra with comultiplication $\Delta:=\unitization  (\hat\sigma\circ\Lambda)$ and $\unitization \exp_\otimes(\psi)=\exp_{(\mathds 1 B,\Delta)}(\hat\psi)$, where $\hat\psi$ is the extension of $\psi$ to a generating functional on $\unitization  B$ (i.e., $\hat\psi(1)=0$). Therefore, $\exp_\otimes\psi$ is a restricted state by the Schoenberg correspondence for $*$-bialgebras \cite{Sch85}.   
\end{proof}

\begin{problem}
  This result provokes the following questions, to which we do not know the answer yet. 
  \begin{enumerate}
  \item Are the non-positive universal products derived from c-bifree-independence half positive or Schoenberg?
  \item If the answer to question (1) is negative, do the reverse implications in Theorem~\ref{main-theorem} hold?
  \end{enumerate}
\end{problem}

\section{Outlook: Brownian motions and interacting Fock spaces for two-faced independences}
\label{sec:gener-brown-moti}

This final section aims to indicate a potential application of the Schoenberg correspondence, without going into much depth. We will restrict ourselves to the case of two-faced universal products (i.e., $2$-$1$-universal products) in order to spare some indices and increase readibility.

\begin{definition}
  Consider a four dimensional vector space $V=\operatorname{span}(x,x^*,y,y^*)$ and the $2$-faced tensor algebra $B=T(V)$ with faces $B^1=T(\mathbb Cx+\mathbb Cx^*)$ and $B^2=T(\mathbb Cy+\mathbb Cy^*)$. This becomes a 2-faced $*$-dual semigroup with the primitive comultiplication $\Lambda(x)=x_1+x_2, \Lambda(y):=y_1+y_2$ and the obvious involution.  The \emph{standard Gaussian cumulant functional}\footnote{The ``Gaussian cumulant functional'' in Example~\ref{ex:Lachs-functor} is closely related to this, but the structure becomes a bit more interesting by allowing non-selfadjoint indeterminates. The variable $x$ in Example~\ref{ex:Lachs-functor} corresponds to $x+x^*$ here.} is the unique linear functional $\psi_{\mathrm g}\colon B\to\mathbb C$  such that $\psi(xx^*)=\psi_{\mathrm g}(xy^*)=\psi(yx^*)=\psi(yy^*)=1$ and $\psi$ vanishes on all other monomials.
\end{definition}

Note that $\psi_{\mathrm g}$ is a restricted generating functional. Therefore, by Schoenberg correspondence, for every positive $2$-faced universal product there is a convolution semigroup of restricted states
\[\varphi^{\odot}_{t}:=\exp_\odot(t\psi_{\mathrm g}).\]

By a general construction (cf.\ \cite[Theorem 4.20 and Section 5.5]{GLS22}), one can associate with every such convolution semigroup of states (uniquely
up to equivalence) a noncommutative probability space $(\mathcal A,\Phi)$ and a non-commutative \emph{L{\'e}vy process}
$j^{\odot}_{s,t}\colon B\to \mathcal A$, $s\geq t\geq0$, with marginal distributions $\Phi\circ j^{\odot}_{s,t}=\varphi^\odot_{s-t}$.

\begin{definition}
  The L{\'e}vy process $(j^{\odot}_{s,t})_{s\geq t}$ with restricted generating functional $\psi_{\mathrm g}$ is called \emph{$\odot$-Brownian motion}.
\end{definition}

For the single-faced independences, these generalized Brownian motions are well-studied, and they can be realized as creation and annihilation processes on the associated Fock spaces (see \cite{Fra06} and references therein).

In the two-faced situation, in particular for $\odot$ one of the ``exceptional'' two-faced universal products from \cite{GVa24}, our general Schoenberg-correspondence promises to be very useful. Those are symmetric universal products which fulfill all the necessary conditions we know of for being positive, but their positivity has not yet been established, because they do not come from universal lifts to the tensor product or the free product in the sense of \cite{GHU23}. Var\v{s}o, in his PhD thesis \cite[Section~6.4]{Var21}, calculated the first few moments $\varphi_1^{\odot}\left((x+x^*+y+y^*)^n\right)$ and showed that the associated Hankel-determinants are positive,\footnote{Actually, Var\v{s}o formulated his statements for selfadjoint indeterminates as in Example~\ref{ex:Lachs-functor}.} which by Schoenberg correspondence can be interpreted as further evidence that the exceptional universal products should be positive. Now, if we assume positivity, by GNS-construction one should be able to define \emph{interacting Fock spaces} (in the sense of \cite[Drf.~18.1]{ALV97} or, more directly, abstract interacting Fock spaces in the sense of \cite[Def.~3.1]{GSk20}) with $j_{s,t}^{\odot}(x^*)$ and $j_{s,t}^{\odot}(y^*)$ as creation operators for the exceptional universal products.

Furthermore, restriction of the two-faced universal product $\odot$ to each of the faces yields two single-faced universal products $\odot_1$ and $\odot_2$, each of which must be one of Muraki's five. It follows that restriction of $(j^{\odot})_{s,t}$ to $B^i$ is equivalent to $(j^{\odot_i}_{s,t})$ and the restrictions to the $*$-algebra generated by $\lambda x + (1-\lambda)y$ for $\lambda\in[0,1]$ gives interpolating generalized Brownian motions. The first example of such interpolating Brownian motions between free and tensor independence was given by Bo\.zejko and Speicher \cite{BoSp92} and received a lot of attention afterwards. Assuming that the free-tensor independence (one of the exceptional independences from \cite{GVa24}) is positive, Schoenberg correspondence gives us another such interpolation for free. 

\section*{Acknowledgements}
I thank Michael Sch{\"u}rmann, Michael Skeide, and Philipp Var\v{s}o for very helpful discussions and feedback on earlier drafts of these notes. I thank Reviewer 2 at \emph{Communications in Mathematical Physics} for pointing at a number of places where giving some additional details hopefully helps readers with a more probabilistic or analytic background to follow the algebraic and categorical arguments. I thank the reviewer at \emph{Journal of Functional analysis} for their thorough reading and in particular for motivating me to write the additional outlook section.  

This work was supported by the German Research Foundation (DFG) grant nos.\ 465189426 and 397960675; it was carried out as a postdoctoral researcher at Saarland University, during the tenure of an ERCIM ``Alain Bensoussan'' Fellowship Programme at NTNU Trondheim, and as a postdoctoral scientific employee at University of Greifswald.

\setlength{\parindent}{0pt}

\end{document}